# Unified Method for Solving General Polynomial Equations of Degree Less Than Five

By Dr. Raghavendra G. Kulkarni

## Abstract

In this paper we present a unified method for solving general polynomial equations of degree less than five.

## 1. Introduction

Is there some common or unified approach for solving general polynomial equations of degree four or less instead of solving them by different methods? The answer is 'yes.' This paper presents a novel method of decomposing the general polynomial equations (of degree four or less) into two constituent lower-degree polynomials as factors, eventually leading to solution. The author uses this method to solve the general quadratic, cubic, and quartic equations. For polynomial equations of degree five and higher (since no general solution is possible in radicals, as shown by Abel and Galois), the method can be used to solve many types of solvable equations [1]. However, since the scope of this paper is limited to general polynomial equations of degree four or less, the solution of polynomial equations of degree five and greater is not discussed here. The following sections contain a brief historical account on polynomial equations, a description of the unified method, and the application of this method to solve quadratics, cubics, and quartics.

## 2. A Brief Historical Account

The quest to solve polynomial equations is not new, and we notice from history that, even before 2000 BC, Greeks, Hindus, and Babylonians knew the solution to quadratics in one form or another. However, the solution of cubics eluded mathematicians for more than three thousand years after the quadratic was solved. In the $11^{\text{th}}$ century, Omer Khayyam solved the cubic geometrically, by intersecting parabolas and circles, and it was Scipione del Ferro, who found the solution to the cubic in 1515, but did not bother to publish it. In 1535, Tartaglia obtained formulas to solve certain types of cubic equations, while Cardano (in 1539) published his solution to the general cubic in his book, *The Great Art* (or *The Rules of Algebra),* using complex numbers, at a time when the use of complex numbers was considered to be absurd. At the same time, Cardano's friend Ferrari published the solution to the general quartic equation. This period witnessed lot of disputes among Scipione del Ferro, Tartaglia, Cardano, and Ferrari, regarding credit for the general solution of the cubic. Several efforts by eminent mathematicians to solve the general quintic on similar lines as adopted for cubics and quartics ended up in failure. In 1770, Lagrange showed that polynomials of degree five or more couldn't be solved by these methods. The same idea, in a more general form, was given by Ruffini in 1799 in his book, *General Theory of Equations* [2, 3]. Further work on the polynomials of degree five and higher was reported by Abel and Galois. In 1826, Abel published the proof of the impossibility of generally solving algebraic equations of a degree higher than four by radicals; while Galois (1832), whose work constituted the foundation of group theory, used this theory to show the same impossibility [3, 4]. This does not mean that there is no algebraic solution to general polynomial equations of degree five and higher. A key contribution towards the algebraic solution of the general quintic came from Bring in 1786, who showed that every quintic can be transformed into a principal quintic of the form: $x^5 + ax + b = 0$. Algebraic solutions to the general quintic and sextic equations have been obtained using the symbolic coefficients [5, 6].

[1]

Coming back to general polynomial equations of degree less than five, we observe that the cubic equation is solved by the well-known Cardano's method, while the quartic equation is solved using Ferrari's method. However, these two methods adopt totally different approaches [2]. Underwood describes a method simpler than Cardano's method for solving the general cubic equation with real coefficients, wherein the cubic polynomial is factored into three factors, and all three roots are computed with comparatively less complexity [7]. There are alternative methods to Ferrari's method to solve the quartic equation [8]. However we notice that there is no *common* approach among the methods used for solving cubics and quartics.

The aim of this paper is to present a *unified* method, which can be adopted to solve quadratics, cubics, and quartics. The proposed method consists of decomposing the given polynomial equation into two constituent lower-degree polynomials, in a novel fashion as explained below.

## 3. The Unified Method

Let the general polynomial equation of degree, $N$, for which solution is sought, be expressed as

$$(1) \qquad x^N + a_{N-1}x^{N-1} + \ldots + a_j x^j + \ldots + a_1 x + a_0 = 0,$$

where $a_j$ ($j = 0$ to $(N-1)$) are real coefficients in the polynomial equation (1) and $N$ is less than five. Our aim is to decompose equation (1) into two factors by means of a common or unified method, irrespective of the degree of equation (1). In order to achieve this objective, the polynomial equation (1) has to be converted into another polynomial equation of the same degree, which is in a convenient form to decompose. For this purpose let us construct a polynomial equation of $N^{th}$ degree as follows:

$$(2) \qquad \frac{[V_M(x)]^K - p^K [W_M(x)]^K}{1 - p^K} = 0,$$

where $V_M(x)$ and $W_M(x)$ are constituent polynomials of degree $M$, where $M < N$, and $p$ is an unknown to be determined. The integer $K$ has to satisfy the relation

$$(3) \qquad\qquad\qquad KM = N,$$

so that polynomial equation (2) is of degree $N$. Also, notice that the integer $K$ must always be greater than unity to facilitate decomposition of equation (2). The polynomials $V_M(x)$ and $W_M(x)$ are given by

$$V_M(x) = x^M + b_{M-1}x^{M-1} + \ldots + b_j x^j + \ldots + b_1 x + b_0$$

$$(4) \quad W_M(x) = x^M + c_{M-1}x^{M-1} + \ldots + c_j x^j + \ldots + c_1 x + c_0,$$

where $b_j$ and $c_j$ ($j = 0$ to $(M-1)$) are coefficients in the polynomials, $V_M(x)$ and $W_M(x)$, respectively. These coefficients ($b_j$ and $c_j$) are unknowns to be determined. Notice that the constructed polynomial equation (2) can be easily decomposed into two factors, and one of the factors is given by

$$(5) \qquad\qquad Y_M(x) = \frac{V_M(x) - pW_M(x)}{1 - p}.$$

Thus if we are able to represent our given polynomial equation (1) in the form of (2), then we can decompose (1) into two factors, and one of the two factors is given by (5). For the polynomial equation (1) to be represented as (2), the corresponding $N$ coefficients in equations (1) and (2) are to be equal. However, notice that the coefficients in equation (2) are not expressed explicitly. Therefore, to bring out the coefficients explicitly, equation (2) is expanded and rearranged in descending powers of $x$ (in the same manner that equation (1) is written), as shown below:

$$x^N + d_{N-1}x^{N-1} + \ldots + d_j x^j + \ldots + d_1 x + d_0 = 0, \tag{6}$$

where $d_j$ ($j = 0$ to $(N-1)$) are coefficients obtained after the expansion and rearrangement of equation (2), and are functions of unknowns, $b_j$, $c_j$, and $p$. Now we are in a position to equate the corresponding coefficients, of the given polynomial equation (1) and the constructed polynomial equation, expressed in the form of (6). This process results in $N$ equations in $(2M + 1)$ unknowns, namely $b_0, b_1, \ldots, b_{M-1}, c_0, c_1, \ldots, c_{M-1}$, and $p$, as shown below:

$$d_j = a_j, \tag{7}$$

where $j = 0$ to $(N-1)$.

If the number of unknowns, $(2M + 1)$, is less than the number of equations, $N$, then some of the coefficients of the given polynomial equation (1) get related to each other, leading to a solution which is not the general solution. To ensure that the proposed method yields the general solution, the number of unknowns $(2M+1)$ has to be equal to the number of equations $N$ when $N$ is odd, and equal to $N + 1$ when $N$ is even, as shown below[1]:

$$2M + 1 = N, \text{ for } N \text{ odd}$$

$$2M + 1 = N + 1, \text{ for } N \text{ even.} \tag{8}$$

The integer $M$ is determined first using (8), and then $K$ is determined from (3) using the values of $N$ and $M$. Observe that when $N$ is even, the number of unknowns is one more than the number of equations. In such case, the extra (one) unknown has to be assigned some convenient value, so that all the unknowns can be determined by solving the $N$ equations given in (7).

---

[1] The reader will note that having at least as many unknowns as equations is a necessary (but certainly not sufficient) condition for the solution of the general $N^{th}$ degree polynomial equation. For example, for $N = 6$ one computes $M = 3$ and $K = 2$, and so $2M + 1 \geq 6$, and yet there exist non-solvable polynomials of degree 6.

With determination of all unknowns, we are in a position to represent the given polynomial equation (1) as a constructed polynomial equation written in the form of (2). Hence the given polynomial equation can be decomposed into two factors, with one of the factors being $Y_M(x)$ as shown in (5). This factor, $Y_M(x)$, when equated to zero results in a polynomial equation of degree $M$ as shown below:

$$(9) \quad x^M + \ldots + \frac{b_j - pc_j}{1-p}x^j + \ldots + \frac{b_1 - pc_1}{1-p}x + \frac{b_0 - pc_0}{1-p} = 0.$$

Solving the $M^{th}$ degree polynomial equation above, we obtain $M$ roots of the given polynomial equation (1). In the same manner, the remaining factor of (2) is equated to zero to obtain a polynomial equation of degree $(N-M)$. By solving this polynomial equation, we obtain the remaining $(N-M)$ roots. Thus all the $N$ roots are determined by factorization of the given polynomial equation into $N$ linear factors, which means that the given polynomial is evaluated over the complex field, **C**. We now demonstrate that this method can be used to solve quadratic, cubic, as well as quartic equations, which are of general nature.

## 4. The Quadratic Equation

We shall now apply the proposed unified method to obtain the roots of the quadratic equation. Let the quadratic equation for which solution is sought be

$$(10) \quad x^2 + a_1 x + a_0 = 0,$$

where $a_0$ and $a_1$ are real, independent, and non-zero coefficients. (When $a_0$ or $a_1$ are zero, the quadratic can readily be factored, and therefore there is no need to proceed further with the proposed method). Applying the unified method to the above quadratic equation $(N=2)$, we obtain the values of $M$ and $K$ from expressions (8) and (3) respectively as $M=1$, and $K=2$. The constructed quadratic equation, as shown in (2), is then given by

$$(11) \quad \frac{[V_1(x)]^2 - p^2 [W_1(x)]^2}{1-p^2} = 0,$$

where $V_1(x)$ and $W_1(x)$ are expressed below (using (4)):

$$V_1(x) = x + b_0$$

$$(12) \quad W_1(x) = x + c_0,$$

and $b_0$ and $c_0$ are unknown coefficients in polynomials, $V_1(x)$ and $W_1(x)$, respectively. Substituting the expressions for $V_1(x)$ and $W_1(x)$ (as given in (12)) into equation (11), we obtain the constructed quadratic equation in the following form:

$$(13) \qquad \frac{(x+b_0)^2 - p^2(x+c_0)^2}{1-p^2} = 0.$$

The equation above (13) is expanded and rearranged in descending powers of $x$ as shown below:

$$(14) \qquad x^2 + \frac{2(b_0 - c_0 p^2)}{1-p^2}x + \frac{b_0^2 - c_0^2 p^2}{1-p^2} = 0.$$

Equating the coefficients of the given quadratic equation (10) with that of the constructed quadratic equation (14), yields the following two equations:

$$(15) \qquad \frac{2(b_0 - c_0 p^2)}{1-p^2} = a_1$$

$$(16) \qquad \frac{b_0^2 - c_0^2 p^2}{1-p^2} = a_0.$$

Notice that there are three unknowns ($b_0$, $c_0$, and $p$) to be determined, while there are only two equations ((15) and (16)) to solve. Therefore one unknown has to be assigned some convenient value. Let us make $c_0 = 0$, and solve equations, (15) and (16), to evaluate $b_0$ and $p$ in terms of $a_0$ and $a_1$. The values of $b_0$ and $p$ obtained are given below:

$$b_0 = \frac{2a_0}{a_1}$$

$$(17) \qquad p = \frac{\left(a_1^2 - 4a_0\right)^{\frac{1}{2}}}{a_1}.$$

Having determined $b_0$ and $p$, the given quadratic equation (10) can be represented as the constructed quadratic equation (13), and hence it can be split into two factors as shown below (noting that $c_0 = 0$):

$$(18) \qquad \frac{x + b_0 - px}{1-p} \cdot \frac{x + b_0 + px}{1+p} = 0.$$

Setting each one of the factors equal to zero in the equation above, we obtain following two linear equations in $x$ as given below:

$$x + \frac{b_0}{1-p} = 0$$

$$x + \frac{b_0}{1+p} = 0.$$

Substituting the values of $b_0$ and $p$ (from (17)) into the two linear equations above and rationalizing the denominators, the two roots, $x_1$ and $x_2$, of the given quadratic equation (10) are obtained as

$$x_1 = \frac{-a_1 - \left(a_1^2 - 4a_0\right)^{\frac{1}{2}}}{2}$$

(19)
$$x_2 = \frac{-a_1 + \left(a_1^2 - 4a_0\right)^{\frac{1}{2}}}{2},$$

which are well-known expressions for the roots of the general quadratic equation. We shall now proceed to solve the cubic equation using the same method.

## 5. The Cubic Equation

Consider the following cubic equation in $x$, for which solution through the proposed unified method is sought:

(20) $$x^3 + a_2 x^2 + a_1 x + a_0 = 0,$$

where $a_0$, $a_1$, and $a_2$ are real and independent coefficients. (Without loss of generality, we assume that $a_2 \neq 0$. If $a_2 = 0$, then replacing $x$ in (20) with $x + r$, $r \neq 0$, will yield a cubic whose $a_2$ term is non-zero, and whose splitting field is identical to the splitting field of (20).) We note that for the case of the cubic equation, $N = 3$, and the integers $M$ and $K$ are determined using expressions (8) and (3) respectively as: $M = 1$, and, $K = 3$. The constructed cubic equation as expressed in (2) is given by

(21) $$\frac{[V_1(x)]^3 - p^3 [W_1(x)]^3}{1 - p^3} = 0,$$

where the polynomials $V_1(x)$ and $W_1(x)$, are given by

$$V_1(x) = x + b_0$$

(22) $$W_1(x) = x + c_0,$$

where $b_0$, $c_0$, and $p$ are unknowns to be determined. Using equations (21) and (22), the constructed cubic equation is further expressed as

$$(23) \qquad \frac{(x+b_0)^3 - p^3(x+c_0)^3}{(1-p^3)} = 0.$$

The equation above is expanded and rearranged in descending powers of $x$, as shown below:

$$(24) \quad x^3 + \frac{3(b_0 - c_0 p^3)}{1-p^3} x^2 + \frac{3(b_0^2 - c_0^2 p^3)}{1-p^3} x + \frac{(b_0^3 - c_0^3 p^3)}{1-p^3} = 0.$$

As described in the unified method, the given cubic equation (20) has to be represented in the form of constructed cubic equation (23) to enable us to decompose equation (20). To accomplish this, we equate the corresponding coefficients in equations, (20) and (24) (since (24) is the expanded form of (23)), resulting in the following three equations in three unknowns ($b_0$, $c_0$, and $p$) as shown:

$$(25) \qquad \frac{3(b_0 - c_0 p^3)}{1-p^3} = a_2$$

$$(26) \qquad \frac{3(b_0^2 - c_0^2 p^3)}{1-p^3} = a_1$$

$$(27) \qquad \frac{b_0^3 - c_0^3 p^3}{1-p^3} = a_0.$$

To determine the unknowns we proceed as follows. Equations (26) and (27) are divided by equation (25) to get the following expressions:

$$(28) \qquad \frac{b_0^2 - c_0^2 p^3}{b_0 - c_0 p^3} = \frac{a_1}{a_2}$$

$$(29) \qquad \frac{b_0^3 - c_0^3 p^3}{b_0 - c_0 p^3} = \frac{3a_0}{a_2}.$$

Equation (25) is rearranged to obtain an expression for $p^3$ as shown:

$$(30) \qquad p^3 = \frac{a_2 - 3b_0}{a_2 - 3c_0}.$$

Using (30) we eliminate $p^3$ from equation (28), and after some algebraic manipulations we get following expression in $b_0$ and $c_0$:

(31) $$(b_0 - c_0)\left[a_1 + 3b_0c_0 - a_2(b_0 + c_0)\right] = 0.$$

The expression above contains two factors. When the factor $(b_0 - c_0)$ is equated to zero, $b_0$ equals $c_0$, and this condition transforms the cubic equation (23) into the cube of a linear polynomial as follows:

$$(x + c_0)^3 = 0.$$

When the cubic above is equated with the given cubic (20), the coefficients, $a_0$, $a_1$, and $a_2$, become interdependent, and we end up with a special case of cubic equation, which is the cube of a linear polynomial. But since we are seeking the solution of the *general* cubic equation, $b_0$ need not be equal to $c_0$ in general. Hence, we have to consider the other factor in (31); and after equating this factor to zero, we get the following relation between $b_0$ and $c_0$ as shown:

(32) $$a_1 + 3b_0c_0 - a_2(b_0 + c_0) = 0.$$

We now eliminate $p^3$ from equation (29), using equation (30), and after some algebraic manipulations, we get the following expression, where again the term $(b_0 - c_0)$ emerges as a factor:

(33) $$(b_0 - c_0)\left[3a_0 + 3b_0c_0(b_0 + c_0) - a_2(b_0^2 + b_0c_0 + c_0^2)\right] = 0.$$

As in the earlier case, we consider the factor other than $(b_0 - c_0)$ in the expression (33), and equate it to zero, and after some manipulations, we obtain following relation between $b_0$ and $c_0$ as shown:

(34) $$3a_0 + a_2 b_0 c_0 + 3b_0 c_0 (b_0 + c_0) - a_2(b_0 + c_0)^2 = 0.$$

It is interesting to see that the unknowns, $b_0$ and $c_0$, appear together as sum and product terms in equations (32) and (34), prompting us to define $f_1$ and $f_2$ as follows:

(35) $$f_1 = b_0 + c_0$$

(36) $$f_2 = b_0 c_0.$$

Re-expressing equations (32) and (34) in terms of $f_1$ and $f_2$, we have

(37) $$a_1 + 3f_2 - a_2 f_1 = 0$$

$$3a_0 + a_2 f_2 + 3 f_1 f_2 - a_2 f_1^2 = 0. \tag{38}$$

By solving equations (37) and (38), $f_1$ and $f_2$ are determined as shown below:

$$f_1 = \frac{a_1 a_2 - 9 a_0}{a_2^2 - 3 a_1} \tag{39}$$

$$f_2 = \frac{a_1^2 - 3 a_0 a_2}{a_2^2 - 3 a_1}. \tag{40}$$

In case the term $\left(a_2^2 - 3 a_1\right)$, which is in the denominator in (39) and (40), is zero, then the given cubic equation (20) is reduced to the special case as shown below:

$$\left(x + \frac{a_2}{3}\right)^3 = (\frac{a_2}{3})^3 - a_0,$$

whose solution can be readily found. Notice from equations, (39) and (40), that the sum $f_1$ and product $f_2$ of $b_0$ and $c_0$ are determined. Therefore, we recognize that $b_0$ and $c_0$ are the roots of the quadratic equation

$$y^2 - f_1 y + f_2 = 0.$$

Thus $b_0$ and $c_0$ are determined as

$$b_0 = \frac{f_1 + \left(f_1^2 - 4 f_2\right)^{\frac{1}{2}}}{2} \tag{41}$$

$$c_0 = \frac{f_1 - \left(f_1^2 - 4 f_2\right)^{\frac{1}{2}}}{2}. \tag{42}$$

The only remaining unknown, $p$, is determined from (30) as shown below:

$$p = \left[\frac{2 a_2 - 3 f_1 - 3 \left(f_1^2 - 4 f_2\right)^{\frac{1}{2}}}{2 a_2 - 3 f_1 + 3 \left(f_1^2 - 4 f_2\right)^{\frac{1}{2}}}\right]^{\frac{1}{3}}. \tag{43}$$

Notice that the constructed cubic equation (23) can be factored as shown below:

(44)
$$\frac{(x+b_0) - p(x+c_0)}{1-p} \cdot \frac{(x+b_0)^2 + p(x+b_0)(x+c_0) + p^2 (x+c_0)^2}{1+p+p^2} = 0.$$

Since we are able to represent the given cubic equation (20) in the form of (23) by determining all the unknowns, it is clear that we have successfully decomposed (20) as shown in (44). Equating each of the factors in (44) to zero, we get two equations as shown below:

$$(45) \qquad (x + b_0) - p(x + c_0) = 0$$

$$(46) \qquad (x + b_0)^2 + p(x + b_0)(x + c_0) + p^2(x + c_0)^2 = 0.$$

Equation (45) is a linear equation in $x$, from which we get one root of the given cubic equation; while equation (46) is a quadratic equation, which provides the other two roots. Thus, solving (45), one root $x_1$ of the given cubic equation is obtained as

$$(47) \qquad x_1 = \frac{pc_0 - b_0}{1 - p}.$$

The remaining two roots, $x_2$ and $x_3$, of the given cubic are obtained when the quadratic equation (46) is solved. The roots are given by

$$(48) \quad \begin{aligned} x_2 &= \frac{-(2b_0 + 2p^2 c_0 + pf_1)}{2 + 2p + 2p^2} \\ &+ \frac{\left[(2b_0 + 2p^2 c_0 + pf_1)^2 - 4(1 + p + p^2)(b_0^2 + p^2 c_0^2 + pf_2)\right]^{\frac{1}{2}}}{2 + 2p + 2p^2} \end{aligned}$$

$$(49) \quad \begin{aligned} x_3 &= \frac{-(2b_0 + 2p^2 c_0 + pf_1)}{2 + 2p + 2p^2} \\ &- \frac{\left[(2b_0 + 2p^2 c_0 + pf_1)^2 - 4(1 + p + p^2)(b_0^2 + p^2 c_0^2 + pf_2)\right]^{\frac{1}{2}}}{2 + 2p + 2p^2} \end{aligned}.$$

We have determined all three roots of the general cubic equation using the unified method. In the last section of the paper a numerical example of a cubic equation is solved, confirming the validity of the formulas obtained for the roots of the cubic equation.

## 6. The Quartic Equation

Consider the following fourth-degree polynomial equation in $x$ as given below:

$$(50) \qquad x^4 + a_3 x^3 + a_2 x^2 + a_1 x + a_0 = 0,$$

where $a_0$, $a_1$, $a_2$, and $a_3$ are real and independent coefficients. To decompose the quartic equation above, using the proposed unified method, we begin by determining the integers $M$ and $K$ using equations (8) and (3), respectively, as $M = 2$ and $K = 2$. The next step is to construct a quartic equation in the form shown in (2), for the case of $N = 4$, $M = 2$, and $K = 2$, as indicated below:

$$\frac{[V_2(x)]^2 - p^2 [W_2(x)]^2}{1 - p^2} = 0, \tag{51}$$

where the constituent polynomials, $V_2(x)$ and $W_2(x)$, are given by:

$$V_2(x) = x^2 + b_1 x + b_0$$

$$W_2(x) = x^2 + c_1 x + c_0, \tag{52}$$

where $b_0$, $b_1$ and $c_0$, $c_1$ are the coefficients of polynomials, $V_2(x)$ and $W_2(x)$, respectively. Substituting the expressions for $V_2(x)$ and $W_2(x)$ from (52) into equation (51), we obtain the constructed quartic equation as shown below:

$$\frac{\left(x^2 + b_1 x + b_0\right)^2 - p^2 \left(x^2 + c_1 x + c_0\right)^2}{1 - p^2} = 0. \tag{53}$$

Expanding and rearranging equation (53) in descending powers of $x$ yields the equation below:.

$$x^4 + \frac{2(b_1 - c_1 p^2)}{1 - p^2} x^3 + \frac{b_1^2 + 2b_0 - (c_1^2 + 2c_0) p^2}{1 - p^2} x^2 + \frac{2(b_0 b_1 - c_0 c_1 p^2)}{1 - p^2} x \tag{54}$$
$$+ \frac{b_0^2 - c_0^2 p^2}{1 - p^2} = 0.$$

Equating the corresponding coefficients of the given quartic equation (50) and the constructed quartic equation (54) results in the four equations given below:

$$\frac{2(b_1 - c_1 p^2)}{1 - p^2} = a_3 \tag{55}$$

$$\frac{b_1^2 + 2b_0 - (c_1^2 + 2c_0) p^2}{1 - p^2} = a_2 \tag{56}$$

$$\frac{2(b_0 b_1 - c_0 c_1 p^2)}{1 - p^2} = a_1 \tag{57}$$

$$\frac{b_0^2 - c_0^2 p^2}{1 - p^2} = a_0. \tag{58}$$

Since there are five unknowns ($b_0$, $b_1$, $c_0$, $c_1$, and $p$) but only four equations to solve, we assign some convenient value to one unknown, so that all the unknowns can be determined. Let us make $c_1 = 0$, and use it in equations that contain $c_1$. As a result, the equations, (55), (56), and (57), are modified as shown below:

$$(59) \qquad \frac{2b_1}{1-p^2} = a_3$$

$$(60) \qquad \frac{b_1^2 + 2b_0 - 2c_0 p^2}{1-p^2} = a_2$$

$$(61) \qquad \frac{2b_0 b_1}{1-p^2} = a_1.$$

(Without loss of generality, we assume that $a_3 \neq 0$. Otherwise, $b_1 = 0$, and consequently, $a_1 = 0$. In this case, (50) is a quadratic in $x^2$, and can be solved as in section 4.) Now we start the process of determining the four unknowns, $b_0$, $b_1$, $c_0$, and $p$, using the four equations, (58), (59), (60), and (61). Dividing the equations, (61), (58), and (60), by equation (59), we obtain following expressions, respectively:

$$(62) \qquad b_0 = \frac{a_1}{a_3}$$

$$(63) \qquad \frac{b_0^2 - c_0^2 p^2}{2b_1} = \frac{a_0}{a_3}$$

$$(64) \qquad \frac{b_1^2 + 2b_0 - 2c_0 p^2}{2b_1} = \frac{a_2}{a_3}.$$

The equation (59) is rearranged to get following expression for $p^2$:

$$(65) \qquad p^2 = \frac{a_3 - 2b_1}{a_3}.$$

From (62) we note that $b_0$ has been evaluated, and from (65) we note that $p^2$ is expressed in terms of $b_1$. Eliminating $p^2$ from equations, (63) and (64), and using equation, (65), we obtain the following two expressions in $b_1$ and $c_0$:

$$(66) \qquad a_3 c_0^2 (2b_1 - a_3) = 2a_0 a_3 b_1 - a_1^2$$

$$(67) \qquad 2c_0 (2b_1 - a_3) = 2a_2 b_1 - 2a_1 - a_3 b_1^2.$$

Observe that both the equations, (66) and (67), contain the factor, $c_0 (2b_1 - a_3)$, on the left hand side. Therefore dividing equation (66) by (67) results in the elimination of factor, $c_0 (2b_1 - a_3)$, and after some rearrangement we get an expression for $c_0$ in terms of $b_1$ as shown below:

$$(68) \quad c_0 = \frac{2\left(2a_0 a_3 b_1 - a_1^2\right)}{a_3\left(2a_2 b_1 - 2a_1 - a_3 b_1^3\right)}.$$

(Notice that in deriving (68), we have tacitly assumed that $c_0\left(2b_1 - a_3\right) \neq 0$. If $c_0\left(2b_1 - a_3\right) = 0$, then either $c_0 = 0$, in which case $W_2\left(x\right)$ factors easily (see (52)), or else $2b_1 - a_3 = 0$, in which case $p = 0$ (see (59)), and hence, (50) can be re-written as $\left[V_2\left(x\right)\right]^2 = 0$ (see 51) and solved using the approach of section 4.) We eliminate $c_0$ from equations, (67) and (68), to obtain the following expression in terms of $b_1$ only:

$$(69) \quad \begin{aligned} a_3^3 b_1^4 &- 4a_2 a_3^2 b_1^3 + 4\left(a_1 a_3^2 + a_2^2 a_3 - 4a_0 a_3\right) b_1^2 \\ &+ 8\left(a_0 a_3^2 + a_1^2 - a_1 a_2 a_3\right) b_1 = 0. \end{aligned}$$

Note that in equation (69) $b_1$ emerges as a factor; however, $b_1 = 0$ is not part of the solution to our general quartic equation, since it makes the value of $p^2$ equal to unity (see (65)), and as a result there is division by zero in equations (53) to (58). Therefore factoring out $b_1$ from equation (69) results in a cubic equation in $b_1$ as shown below:

$$(70)$$
$$b_1^3 - \frac{4a_2}{a_3} b_1^2 + \frac{4\left(a_1 a_3 + a_2^2 - 4a_0\right)}{a_3^2} b_1 + \frac{8\left(a_0 a_3^2 + a_1^2 - a_1 a_2 a_3\right)}{a_3^3} = 0.$$

Solving the cubic equation above, either by Cardano's method or by the unified method described above, $b_1$ is determined. The remaining unknowns, $p^2$ and $c_0$, are determined from equations, (65) and (68) respectively.

Since all unknowns have been determined, we are in a position to represent the given quartic equation (50) in the form of our constructed quartic equation (53). Therefore we can decompose the given quartic equation into two factors as

$$(71)$$
$$\frac{\left(x^2 + b_1 x + b_0\right) - p\left(x^2 + c_0\right)}{1 - p} \cdot \frac{\left(x^2 + b_1 x + b_0\right) + p\left(x^2 + c_0\right)}{1 + p} = 0.$$

Equating each factor in (71) to zero, we obtain two quadratic equations as shown below:

$$(72) \quad x^2 + \frac{b_1}{1 - p} x + \frac{b_0 - c_0 p}{1 - p} = 0$$

$$(73) \qquad x^2 + \frac{b_1}{1+p}x + \frac{b_0 + c_0 p}{1+p} = 0.$$

We determine all four roots ($x_1$, $x_2$, $x_3$, and $x_4$) of the given quartic equation (50), by solving the quadratic equations, (72) and (73). Quadratic equation (72) provides two roots ($x_1$ and $x_2$) as follows:

$$x_1 = \frac{-b_1 + [b_1^2 - 4(b_0 - c_0 p)(1-p)]^{\frac{1}{2}}}{2(1-p))}$$

$$(74) \qquad x_2 = \frac{-b_1 - [b_1^2 - 4(b_0 - c_0 p)(1-p)]^{\frac{1}{2}}}{2(1-p))}.$$

The other two roots, $x_3$ and $x_4$, are obtained by solving the quadratic equation (73) as shown below:

$$x_3 = \frac{-b_1 + [b_1^2 - 4(b_0 + c_0 p)(1+p)]^{\frac{1}{2}}}{2(1+p)}$$

$$(75) \qquad x_4 = \frac{-b_1 - [b_1^2 - 4(b_0 + c_0 p)(1+p)]^{\frac{1}{2}}}{2(1+p)}.$$

We have solved the general quartic equation and determined all four roots using the proposed unified method. In the next section, we offer some numerical examples of cubic and quartic equations which are solved using the unified method.

## 7. Numerical Examples

We consider an example of a cubic equation having irrational and complex roots, and use the unified method to extract all three roots. The cubic equation chosen is

$$x^3 - 2.049888 x^2 + 3.1010205 x + 11.313708 = 0.$$

The parameters, $f_1$ and $f_2$ are determined from the coefficients of the preceding cubic equation as

$$f_1 = 21.20755, \ f_2 = -15.52471,$$

using equations (39) and (40). The values of $b_0$, $c_0$, and $p$ are obtained from (41), (42), and (43) respectively as

$$b_0 = 21.91592, \ c_0 = -0.708376, \ \text{and} \ p = -9.658787.$$

The three roots are then determined from (47), (48), and (49) as

$$x_1 = -1.4142, \ x_2 = 1.73205 + 2.23607i, \ x_3 = 1.73205 - 2.23607i.$$

We now turn our attention to a numerical example of a quartic equation, whose roots are irrational and complex, in order to represent the most general case. The quartic equation chosen is as follows:

$$x^4 + 2.0533927x^3 - 2.8917903x^2 + 7.6758959x + 29.5803989 = 0.$$

The parameters, $b_0$, $b_1$, $c_0$ and $p$ are determined from equations, (62), (70), (68), and (65) respectively, as given below:

$$b_0 = 3.738154, \ b_1 = -13.44884, \ c_0 = 5.336053, \text{ and } p = 3.754882.$$

Using these values, the four roots are evaluated from the expression sets, (74) and (75), as

$$x_1 = -2.236067, \quad x_2 = -2.645752, \ x_3 = 1.414213 + 1.732051i,$$
$$\text{and} \quad x_4 = 1.414213 - 1.732051i.$$

## 8. Conclusions

We have presented a unified method to solve general polynomial equations of degree four or less. The formulas obtained for the roots of cubic and quartic equations have been verified in the preceding examples.

## Acknowledgments

The author thanks the management of Bharat Electronics for supporting this work.

The editor thanks Rob Underwood of Auburn University at Montgomery for his insightful comments.

## Editor's Note

The reader will no doubt wonder what happens when this unified method is applied to a polynomial equation of degree greater than 4. The author answers this question, at least in part, in [1]. There, he applies the unified method to the quintic equation $x^5 + a_4 x^4 + a_3 x^3 + a_2 x^2 + a_1 x + a_0 = 0$. The overall strategy involves converting the quintic into a sextic by multiplying both sides of the equation by the factor $(x + k)$, introducing the root $x = -k$ in the process. In the fashion of equation (11), the resulting sextic

$$\begin{aligned}x^6 \ &+ \ (a_4 + k)\,x^5 + (a_3 + a_4 k)\,x^4 + (a_2 + a_3 k)\,x^3 \\ &+ \ (a_1 + a_2 k)\,x^2 + (a_0 + a_1 k)\,x + a_0 k = 0\end{aligned}$$

is re-written in the form

$$\frac{\left(x^3 + b_2 x^2 + b_1 x_1 + b_0\right)^2 - p^2 \left(x^3 + c_2 x^2 + c_1 x_1 + c_0\right)^2}{1 - p^2} = 0,$$

and then decomposed into the product of cubic factors

$$\frac{\left(x^3 + b_2 x^2 + b_1 x_1 + b_0\right) - p\left(x^3 + c_2 x^2 + c_1 x_1 + c_0\right)^2}{1 - p}$$

$$\cdot \frac{\left(x^3 + b_2 x^2 + b_1 x_1 + b_0\right) + p\left(x^3 + c_2 x^2 + c_1 x_1 + c_0\right)^2}{1 + p} = 0.$$

Comparison of coefficients leads to six equations in the eight unknowns $b_0, b_1, b_2, c_0, c_1, c_2, p$, and $k$. *The choice of auxiliary equations (e.g., $b_2 = 0$ and $c_2 = 0$) determines the type of solvable quintic that can be solved using this technique.* From this point forward, the roots of each cubic are found using the approach in section 5.

Deputy General Manager, HMC division, Bharat Electronics
Jalahalli Post,
Bangalore-560013, INDIA.
Phone: +91-80-22195270, Fax: +91-80-28382738
e-mail: rgkulkarni@ieee.org